\documentstyle[12pt,amssymb]{amsart}
\newmathalphabet*{\frakb}{euf}{b}{n}
\newmathalphabet*{\eusb}{eus}{b}{n}
\size{10}{18pt}\selectfont
\begin{document}
\title{A POSSIBLE ORIGIN OF LOGIC}
\author[V.Bulitko]{V.K.Bulitko\\
\vspace{.2cm}
{\tiny Odessa State University,\\ 
Fontanskaya Doroga 67-40, Odessa 65062, Ukraine\\
e-mail: booly@@te.net.ua}
}
\vspace*{-2.0cm}
\maketitle

\begin{center} {\small Subj-class: Logic (LO)} \end{center}

\section{introduction}

An origin is often an intriguing issue. 
It becomes doubly intriguing when the logical form of thinking
is considered. In this paper we will investigate exactly that: 
we will conjecture on the origin of basic instruments of logical
thinking.
 
As the starting point we consider a very general case 
of a subject interacting 
with an environment. The subject might often have 
a choice of actions to take. 
Additionally, it is of fundamental importance that 
as a rule the subject's internal representations of (beliefs about) the 
environment are inaccurate.

Beginning with this very general schema we will eventually arrive 
at a possible mechanism of
how Logic could have emerged. Our reasoning will be supported with 
some in-depth case studies.

To formalize the settings we introduce two
partially ordered sets $(M,\le_m)$ of environmental positions 
that the subject can 
occupy with its "internal" ordering of their values. A 
partially ordered set $(L,\le_l)$ of "external" estimates 
can be used to characterize the objective value of an environmental position.
For that we have an estimate function 
$\psi:(M,\le_m)\to(L,\le_l)$. So the partial order $\le_m$ 
describes the subject's internal (and therefore subjective)
representations of desirability or reachability
of positions from $M$. 
Then $(\psi(M),\le_l)$ yields the objective ("external")
values of the positions. 

A crucial task of the subject is to find a position $m\in M$ maximizing
the estimate $\psi(m)$ in the poset $(\psi(M),\le_l)$. 

No logical inference is needed 
to find an extremum with respect to its subjective 
("internal") ordering $(M,\le_m)$. It is possible
to reach such a maximum by means of a greedy algorithm. 
The same is true when $\psi$ is a monotonical function 
(i.e., the condition $x\le y\implies\psi(x)\le\psi(y)$ is satisfied).

However, the task becomes more complex when the estimate
function is not monotonical. 
It is then most naturally to explore some version of the successive 
approximation method. 

The rest of the paper is organized as follows. We will start with the general
theory, move onto a special case, and then apply the introduced approach to the 
classical two-valued propositional
logic and modal propositional logics. A conclusion section will follow. We 
will use the standard mathematical logic and 
partial order set theory notation \cite{1,2}.

\section{General theory}

The suggested version operates on a particular yet very general 
representation of the operator $\psi$ in the so-called 
"approximating form". 
It uses three axiomatically defined operations 
$\boxminus,\boxplus,\circledcirc$ based only on general properties
of the posets $(M,\le_m),(L,\le_l)$ as follows.

For every poset $(R,\le_r)$ the standard mappings  
$(\cdot)^{\vartriangle},(\cdot)^{\triangledown}:R\to2^R$ 
are defined by $t^{\vartriangle}=\{t'\in R|t'\le_r t\},
t^{\triangledown}=\{t'\in R|t\le_r t'\}.$

Let us suppose a binary operation  
$\boxminus:L\times L\to L$ and unary operations
$\boxplus:2^L\to L,\circledcirc:L\to L$ are defined in such a way 
that the following system $\cal A$ of axioms holds. 
\begin{description}
\item[$\cal A_1$]  
       $(\forall S\subseteq M)(\exists\breve S\subseteq S)
       [(\forall s\in S)(\exists\breve s\in\breve S)
       [\breve s\le_ms]\ \&\ (\forall\breve s,\breve s'\in\breve S)
       [\breve s\not<_m\breve s']]]$.
\item[$\cal A_2$]  
       $(\forall L',L''\subseteq L)(\forall x\in L')[(x\le_l\boxplus(L'))\&
       (L'\subseteq L''\implies\boxplus(L')\le_l\boxplus(L''))]$.
\item[$\cal A_3$]  
       $(\forall l,l'\in L)[\boxminus(l,\circledcirc(l))=l\ \&\ 
       (l\le_ll'\implies\circledcirc(l)\le_l\circledcirc(l'))]$.
\item[$\cal A_4$]
       $(\forall l,l'\in L)[l\le_ll'\implies(\exists l''\in L)
       [\boxminus(l',l'')=l\ \&\ \circledcirc(l')\le_ll'']]$.
\end{description}

For every operator $\nu:M\to L$ we call set $\frak n(\nu)=
\{(m,m')|(m\le_mm')\ \&\ (\nu(m)\not\le_l\nu(m'))\}$ {\sl non-monotonicity
domain of $\nu$}. 
If $\frak n(\nu)=\emptyset$ then $\nu$ is called monotonical operator.

{\sl{\bf Theorem 1.} Let all axioms of the system $\cal A$ be satisfied
for $(M,\le_m),(L,\le_l)$ and $(M,\le_m)$ have only finite increasing
chains.  Then for every $\psi:M\to L$ there exists a representation
$\psi=\boxminus(\varphi_1,
\boxminus(\varphi_2,\boxminus(\varphi_3,\dots)))$ where all 
$\varphi_i,i=1,2,3\dots$ are monotonical mappings from  
$(M,\le_m)$ to $(L,\le_l)$.

The number of occurences of the operation $\boxminus$ in this 
representation does not exceed the maximal length among lengths
of all increasing chains in poset $(M,\le_m)$.
}

{\sl\bf Proof.} 
Let us reduce the problem for given operator    
$\psi$ to the same problem for an simpler operator
$\psi_1$ such that the following holds
$\psi=\boxminus(\varphi_1,\psi_1)$ and $\frak n(\psi_1)
\subsetneqq\frak n(\psi)$.

First we define $M_1=\{x\in M|\frak n(\psi)\cap(x^{\vartriangle}
\times x^{\vartriangle})\neq\emptyset\}$,
$M^1=\overline{M_1}$ ¨ 
\begin{center}
$\varphi_1(x) =
\begin{cases}
\boxplus(\psi(x^{\vartriangle})),&  x\in M_1,\\
\psi(x),&   x\in M^1.
\end{cases}$
\end{center}
Then we set $\psi_1(x)$ to any such $z\in L$ that
$\boxminus(\varphi_1(x),z)=\psi(x)\ \&\ 
\circledcirc(\varphi_1(x))\le_lz$ if
$\varphi_1(x)\neq\psi(x)$. Otherwise we set 
$\psi_1(x)=\circledcirc(\psi(x))$.

Existence of the element $z$
in the definition is guaranteed by axioms $\cal A_3,\cal A_4$.
Now equality $\psi(x)=\boxminus(\varphi_1(x),\psi_1(x))$ is true
because of the definitions of $\varphi_1,\psi_1$. 

Let us prove that operator $\varphi_1:(M\le_m)\to(L,\le_l)$ is
monotonical one. 

First,    
$\varphi_1=\psi$ on $M^1$ and we may use condition
$x,y\in M^1\& x\le_my\implies\psi(x)\le_l\psi(y)$. 
Indeed, otherwise 
$\psi(x)\not\le_l\psi(y),x\le_my,\psi(x)\neq\psi(y)$ and therefore
$y\in M_1\cap M^1$. However, $M^1\cap M_1=\emptyset$ 
which leads to a contradiction.  

Second, $\varphi_1$ maps $(M_1,\le_m)$ into
$(L,\le_l)$ monotonically in accordance with $\cal A_2$. 

Finally, let us consider the "mixed" case when $x\in M^1,y\in M_1$
and all elements of $M$ are comparable with respect to $\le_m$. 
It is clear $y\le_mx$ is impossible since condition
$z\in M_1\implies z^{\triangledown}\subseteq M_1$ follows 
from the definition of $M_1$ immediately.

Thus, it remains to consider the possibility of $x\le_my$. 
In this case
$\varphi_1(y)=\boxplus(\psi(y^{\vartriangle}))\ge_l\psi(x)$
in accordance to $\cal A_2$. On the other hand, $\psi(x)=\varphi_1(x)$
on $M^1$ follows from the definition of $\varphi_1$. 
Hence operator $\varphi_1$ is monotonic. 

We are now ready to prove the last assertion of the theorem. For that is is sufficient
to show the inclusion $M^1\cup\breve{M_1}\subseteq M^2$.
Here  $M^2,M_2$ 
are defined for $\psi_1$ in the same way as $M^1,M_1$ were defined
for $\psi$ above. $\breve{M_1}$ is the set of all minimal elements 
of set $M_1$, see $\cal A_1$. 
Namely: $M^2=\overline{M_2}$ and 
$M_2=\{x\in M|\frak n(\psi_1)\cap(x^{\vartriangle}
\times x^{\vartriangle})\neq\emptyset\}$. 

From here we have 
$M_2\subseteq(M_1\setminus\breve{M_1})$ and   
$\frak n(\psi_1)\subseteq\frak n(\psi)\setminus\breve M_1\times M_1$. 
So the sequence $M_1\supsetneqq M_2\supsetneqq M_3\supsetneqq\dots$ 
interrupts on a step with the number that can not be higher the highest
of lengths of the increasing chains in poset $(M,\le_m)$. Indeed, since
$\breve M_2\subseteq M_1\setminus\breve M_1$ then in accordance
with $\cal A_1$
for every element $y\in\breve M_2$ there exists some $x\in\breve M_1$
such that $x<_ly$. Therefore, one can choose some increasing chain 
of represetatives of sets $\breve M_1,\breve M_2,\breve M_3,\dots$
which are mutually disjoint sets. 

We will now prove $M^1\cup\breve{M_1}\subseteq M^2$.
First, $\varphi_1(x)=\psi(x)$ is true for every $x\in M^1$.
From here $\psi_1(x)=\circledcirc(\psi(x))$.
However, mapping $\psi_1$ is 
monotonical on $M^1$ in view of $\cal A_3$ and
since $\psi$ is monotonic on $M^1$. 
So $(M^1\times M^1)\cap\frak n(\psi_1)=\emptyset$ and therefore 
$M^1\subseteq M^2$.

Further, let $x,y\in M^1\cup\breve{M_1}$ and $x\le_m y$. Then
we can show that $\psi_1(x)\le_l\psi_1(y)$. 
Indeed, the case $x,y\in M^1$ was considered above. 
The case $x,y\in\breve{M_1}$ is impossible since all elements of 
$\breve{M_1}$ are incomparable by the definition.
We saw above that $x\in M_1\ \&\ x\le_my\implies y\in M_1$. 
Besides $M^1\cap M_1=\emptyset$. 
Therefore, $x\in M^1,y\in\breve M_1$ is the only case remaining 
to consider. 
By definition  $\psi_1(x)=\circledcirc(\psi(x))$ and 
relation $\boxminus(\varphi_1(y),\psi_1(y))=\psi(y)$ holds. Moreover, 
$\psi(y)<_l\varphi(y)$. In accordance with $\cal A_4$ we have 
$\circledcirc(\varphi_1(y))\le_y\psi_1(y)$. Hence
$\psi_1(x)\le_l\psi_1(y)$ takes place since $\circledcirc$ is a monotonical
operation in view of $\cal A_3$ and 
$\psi(z)\le_l\varphi_1(z),z\in M$ in accordance to $\cal A_2$ and the 
construction. \ \ \ $\Box$

Instead of or together with $\cal A$ the dual axiom system
$\cal A^{\star}$ can be fulfilled. It is obtained by replacing
$\le$ with $\ge$ and $\boxplus,\boxminus,\circledcirc$ with
$\boxplus^{\star},\boxminus^{\star},\circledcirc^{\star}$ correspondingly:
\begin{description}
\item[$\cal A^{\star}_1$]  
       $(\forall S\subseteq M)(\exists\breve S\subseteq S)
       [(\forall s\in S)(\exists\breve s\in\breve S)
       [\breve s\ge_ms]\ \&\ (\forall\breve s,\breve s'\in\breve S)
       [\breve s\not<_m\breve s']]]$.
\item[$\cal A^{\star}_2$]  
       $(\forall L',L''\subseteq L)(\forall x\in L')[(x\ge_l
       \boxplus^{\star}(L'))\&
       (L'\subseteq L''\implies\boxplus^{\star}(L')\ge_l
       \boxplus^{\star}(L''))]$.
\item[$\cal A^{\star}_3$]  
       $(\forall l,l'\in L)[\boxminus^{\star}(\circledcirc^{\star}(l),l)=l
       \ \&\ 
       (l\ge_ll'\implies\circledcirc^{\star}(l)\ge_l
       \circledcirc^{\star}(l'))]$.
\item[$\cal A^{\star}_4$]
       $(\forall l,l'\in L)[l\ge_ll'\implies(\exists l''\in L)
       [\boxminus^{\star}(l'',l')=l\ \&\ 
       \circledcirc^{\star}(l')\ge_ll'']]$.
\end{description}
Then the dual theorem holds:

{\sl{\bf Theorem 1$^{\star}$.} Let all axioms of the system $\cal A^{\star}$ 
be fulfilled for posets $(M,\le_m),(L,\le_l)$, 
and operators $\boxplus^{\star},\boxminus^{\star},\circledcirc^{\star}$
and  $(M,\le_m)$ have only finite decreasing chains.  
Then for every operator $\psi:M\to L$ there exists representation
$\psi=\boxminus^{\star}(\dots\boxminus^{\star}(\boxminus^{\star}
(\varphi_{n+1},\varphi_{n}),\linebreak
\varphi_{n-1})\dots,\varphi_1)$ where all
$\varphi_i,i=1,\dots,n,n+1,$ are monotonical mappings from
$(M,\le_m)$ to $(L,\le_l)$.

The number $n$ of occurences of operations $\boxminus^{\star}$ 
in the representation does not exceed 
the highest length among the lengths of
decreasing chains in $(M,\le_m)$.
}

We call the representing forms from these theorems 
{\sl approximating forms}. 
Another way to obtain approximating forms is suggested in theorem 2 below.

Let us suppose a binary operation  
$\boxminus,\uplus:L\times L\to L$ and unary operations
$\circledcirc:L\to L$ are defined in such a way 
that the following system $\cal B$ of axioms takes place. 
\begin{description}
\item[$\cal B_1$]  
       $(\forall S\subseteq M)(\exists\breve S\subseteq S)
       [(\forall s\in S)(\exists\breve s\in\breve S)
       [\breve s\le_ms]\ \&\ (\forall\breve s,\breve s'\in\breve S)
       [\breve s\not<_m\breve s']]]$.
\item[$\cal B_2$]  
       $(\forall x,y\in L)[x,y\le_l\uplus(x,y)]$.
\item[$\cal B_3$]  
       $(\forall l,l'\in L)[\boxminus(l,\circledcirc(l))=l\ \&\ 
       (l\le_ll'\implies\circledcirc(l)\le_l\circledcirc(l'))]$.
\item[$\cal B_4$]
       $(\forall l,l'\in L)[l\le_ll'\implies(\exists l''\in L)
       [\boxminus(l',l'')=l\ \&\ \circledcirc(l')\le_ll'']]$.
\end{description}

{\sl{\bf Theorem 2.} Let all axioms of the system $\cal B$ be satisfied
for $(M,\le_m),(L,\le_l)$ and $(M,\le_m)$ have only finite increasing
chains.  Then for every $\psi:M\to L$ there exists a representation
$\psi=\boxminus(\varphi_1,
\boxminus(\varphi_2,\boxminus(\varphi_3,\dots)))$ where all 
$\varphi_i,i=1,2,3\dots$ are monotonical mappings from  
$(M,\le_m)$ to $(L,\le_l)$.

The number of occurences of the operation $\boxminus$ in this 
representation does not exceed the maximal length among the lengths
of all increasing chains in poset $(M,\le_m)$.
}

{\sl\bf Proof.} First, in the case when 
$(\forall x\in M)[|x^{\vartriangle}|<\infty$ is true
we can prove our theorem using theorem 1. 
For that we only need to note that in this case
it is possible to replace $\boxplus(\psi(x^{\vartriangle})$ with 
any expression of kind $\uplus(psi(z_1),uplus(\dots
\uplus(\psi(z_{n_1},\psi(z_{n})\dots))$. Here $z_1,\dots,z_n$ is an
enumeration of the finite set $x^{\vartriangle}$. 
Indeed, in the proof of theorem 1 we used axiom $\cal A_2$ 
only for subsets of $L$ of the form $\psi(x^{\vartriangle})$. 
Thus, it is sufficient to check only that axiom $\cal A_2$ is
true for sets of kind $\psi(x^{\vartriangle})$. 
This check is a trivial one on the base of axiom $\cal B_2$ 
for operation $\uplus$. 

Otherwise, when there are infinite sets $x^{\vartriangle}$ we can make use
of the condition of fineteness of increasing chaines in $(M,\le_m)$. 
Let us associate every non-minimal element 
$x\in M$ with some maximal with respect to 
the inclusion relation $\subseteq$ increasing 
chain $x_1<_mx_2<_m\dots<_mx_{k+1}=x$. So $x_1$ is the
minimal element of $(M,\le_m)$ and for any $y,j\in\{1,k\}$ if 
$x_j\le_my\le_x{j+1}$ then $x_j=y\vee y=x_{j+1}$. Let us then denote
the previous element $x_k$ of the chain by $\hat x$. 

Now we replace the definition of operator $\varphi_1$ from the proof
of theorem 1 above with the following inductive definition:\\
{\sl Basis:} $x\in M^1$. Then $\varphi_1(x)=\psi(x)$.\\
{\sl Induction Step:} $x\in M_1$ and $\hat x)$ is defined. 
Then we set $\varphi_1(x)=\uplus(\psi(x),\varphi_1(\hat x))$. 

From this we evidently have that $\varphi_1$ is a monotonical operator
and $\psi(x)\le_l\varphi_1(x),x\in M$. The remaining part of the 
proof follows the corresponding part of theorem 1 proof. $\ \ \Box$

The dual theorem relates with the dual axiom system $\cal B^{\star}$.
\begin{description}
\item[$\cal B^{\star}_1$]  
       $(\forall S\subseteq M)(\exists\breve S\subseteq S)
       [(\forall s\in S)(\exists\breve s\in\breve S)
       [\breve s\ge_ms]\ \&\ (\forall\breve s,\breve s'\in\breve S)
       [\breve s\not<_m\breve s']]]$.
\item[$\cal B^{\star}_2$]  
       $(\forall x,y\in L)[x,y\ge_l\uplus^{\star}(x,y)]$.
\item[$\cal B^{\star}_3$]  
       $(\forall l,l'\in L)[\boxminus^{\star}(\circledcirc^{\star}(l),l)=l
       \ \&\ 
       (l\ge_ll'\implies\circledcirc^{\star}(l)\ge_l
       \circledcirc^{\star}(l'))]$.
\item[$\cal B^{\star}_4$]
       $(\forall l,l'\in L)[l\ge_ll'\implies(\exists l''\in L)
       [\boxminus^{\star}(l'',l')=l\ \&\ 
       \circledcirc^{\star}(l')\ge_ll'']]$.
\end{description}
Then the dual theorem holds:

{\sl{\bf Theorem 2$^{\star}$.} Let all axioms of the system $\cal A^{\star}$ 
be fulfilled for posets $(M,\le_m),(L,\le_l)$, 
operators $\uplus^{\star},\boxminus^{\star},\circledcirc^{\star}$
and  $(M,\le_m)$ have only finite decreasing chains.  
Then for every operator $\psi:M\to L$ there exists representation
$\psi=\boxminus^{\star}(\dots\boxminus^{\star}(\boxminus^{\star}
(\varphi_{n+1},\varphi_{n}),\linebreak
\varphi_{n-1})\dots,\varphi_1)$ where all
$\varphi_i,i=1,\dots,n,n+1,$ are monotonical mappings from
$(M,\le_m)$ to $(L,\le_l)$.

The number $n$ of occurences of operations $\boxminus^{\star}$ 
in the representation does not exceed the highest length among the 
lengths of decreasing chains in $(M,\le_m)$.
}

\subsection{A special case}

Sometimes it is possible to choose another type of 
operators $\varphi_i$ in the 
previous theorems. We suggest some condtions for that in the following new 
axiom system $\cal B^+$ which consists of 
the above-introduced system $\cal B$ 
completed with the following axiom:
\begin{description}
\item[$\cal B^+_5$]  
$(\forall l\in L)(\exists\hat l\in L)[l\le_l\hat l\ \&\ \neg(\exists l'\in L)
[\hat l<_ll']].$
\end{description}
This axiom postulates that for any element $l\in L$
there exists at least one maximal element 
$\hat l$ of $L$ greater than $l$. 

Let us denote by 
$\min M$ the class of all minimal elements of $(M,\le_m)$, and by
$\max L$ the class of all maximal elements of $(L,\le_l)$. 
Then let us denote by $\theta:M\to L$ any such function
that for every increasing chain $m_1<_mm_2<_m\dots<_mm_t$ where $m_1\in
\min M$ the following conditions are satisfied. \\
1) $\theta (m_i)=\circledcirc(m_)\vee\theta (m_i)\in\max L\ \&\ 
m_i\le_l\theta (m_i)$;\\
2) $\theta (m_i)\in\max L\ \&\ i\le j\implies \theta (m_j)\in\max L$.
This condition means that $\theta (m)\in\max L\implies(\forall x\in
m^{\triangledown})[\theta(x)\in\max L]$. 
 
We call these mappings $\theta$-{\sl mappings}.

At last, let $\cal K$ be a class consisting of $\theta$-mappings such
that for any pair $(y,x)\in<_m$ there exists a 
$\theta$-mapping $\theta_{y,x}$ obeying the conditions  
$\forall y'[y'\notin x^{\triangledown}\implies\theta_{y,x}(y')=
\circledcirc(y')], \forall x'[x'\in x^{\triangledown}\implies
\theta_{y,x}(x')\in\max L]$. We refer to such functions as  
{\sl special $\theta$-functions}.

{\sl{\bf Theorem 3} Let all axioms of the system $\cal B^+$ be satisfied
for $(M,\le_m),(L,\le_l)$, $\cal K$ satisfy the condition above, 
and $(M,\le_m)$ be finite.
Then for every $\psi:M\to L$ there exists a formula $\Phi(z_1,\dots,z_n)$ 
that
only operations $\boxminus,\uplus$ occur and there exists a substitution
$p:\{z_1,\dots,z_n\}\to\cal K$ such that 
$\psi=Sb^{z_1\ \ \dots\ \ z_n}_{p(z_1)\dots p(z_n)}\Phi(z_1,\dots,z_n)$. 
}

{\bf\sl Proof.} 
We will follow theorem 1 proof but re-define $\varphi_i$. First, we
choose a pair $(y,x)\in<_m$ such that $y$ is a maximal in $M^1$ and $x$ 
immideately follows $y$ in $(M,\le_m)$. Then define $\varphi_{y,x}$ as: 
\begin{center}
$\varphi_{y,x}(z) =
\begin{cases}
\uplus(\psi(z),\theta_{y,x}),&  z\in M_1,\\
\psi(z),&   z\in M^1.
\end{cases}$
\end{center}
Then we set $\psi_{y,x}(u)$ equal to any such $z\in L$ that
$\boxminus(\varphi_{y,x}(u),z)=\psi(u)\ \&\ 
\circledcirc(\varphi_{y,x}(u))\le_lz$ if
$\varphi_{y,x}(u)\neq\psi(u)$. Otherwise we set 
$\psi_{y,x}(u)=\circledcirc(\psi(u))$.

Analogously to the proof of theorem 2 it can be shown that 
$\psi=\boxminus(\varphi_{y,x},\psi_{y,x})$ where $\frak n(\psi_{y,x}),
\frak n(\varphi_{y,x})\subsetneqq\frak n(\psi)$. 

On further steps we handle
$\psi_{y,x},\varphi_{y,x}$ in the same way and so forth. 
Since $M$ is a finite set and we use special $\theta$-functions 
this reduction converges in a finite number of steps. 
$\ \ \ \Box$

Of course the last theorem can be reformulated in the dual form.

\section{Consideration of the classical two-valued propositional logic
from the developed approach}

It is easy to arrive at the classical two-valued propositional logic now. 
For that it is sufficient to choose
$(\{0,1\},0\le1$ as $(L,\le_l)$ and 
the standard poset $(\cal B^n,
\preccurlyeq)$ on boolean cube $\cal B^n$ as poset $(M,\le_m)$.
It is well known that every finite poset can be isotonically 
included into $(\cal B^n,\preccurlyeq)$ for the appropriate $n$. 

Also it is well known that poset $(\cal B^n,\preccurlyeq)$
is a self-dual poset for any  $n$. Therefore, both 
above-introduced representations
take place in this case. 

{\sl{\bf Lemma.} 1) The system of posets $(\cal B^n,\preccurlyeq),
(\cal B,\le)$ as $(M,\le_m),(L\le_l)$ correspondingly and operation
$\to$ as $\boxminus^{\star}$, operation $\Bbb{\bold1}:\cal B^n\to\{1\}$ as
$\circledcirc^{\star}$, and operation  
$\underset{\vec{\beta}\preccurlyeq\vec{\alpha}}{\&}\vec{\alpha}$
as $\boxplus^{\star}(\beta^{\triangledown})$ fulfil axiom
set $\cal A^{\star}$.\\
2) The system of posets $(\cal B^n,\preccurlyeq),
(\cal B,\le)$ as $(M,\le_m),(L\le_l)$ correspondingly and operation
$\to$ as $\boxminus^{\star}$, operation $\Bbb{\bold1}:\cal B^n\to\{1\}$ as
$\circledcirc^{\star}$, and operation $\&$
as $\uplus^{\star}$ fulfil axiom set ${\cal B^+}^{\star}$. 
}

{\sl\bf Proof.} This can done via a routine check of the axioms.

The direct corollary of this lemma and  theorems above is

{\sl{\bf Theorem 4.} In the special cases  
of finite "internal" orders $(M,\le_m)$
and linear "external" orders $(L,\le_l),|L|=2$,  
approximating forms from every of theorems 1,2,3 and their dual ones 
generate all
formulas of the classical propositional logic (within logical equivalence). 
}

Also the following interesting statement follows.

{\sl{\bf Corollary.} Every $n$-argument logical (boolean) function 
$f$ can be represented by the implicative normal form 
$f=P_k\to P_{k-1}\to\dots\to P_1$, where 
$k\le n$, and $P_i,i=\overline{1,k},$ are monotonical boolean function.
}

It is remarkable that just the dual approximating form presents the usual 
propositional implication or that operation $\to^{\star}$ is not
presented in the natural language. In our opinion, the main reason 
is that our dual approximating forms of theorems 1$^{\star}$, 2$^{\star}$ 
start from 
a given operator $\psi$ and approximate it by means of successive 
simplifications: $\psi_1=\boxminus^{\star}(\psi,\varphi_1),\psi_2=
\boxminus^{\star}(\psi_1,\varphi_2),\dots$
while $\psi_i$ is not a monotonical operator (i.e. not an "easy" one). 
Thus, the approximation starts from the target unlike in  
the case of the approximating form from theorems 1,2. 

Now one can consider the classical two-valued
propositional logic merely as a realization of 
the above-mentioned principle of successive approximations for 
the decision-making problems within subject-environment survival framework. 

Thus, from this viewpoint, 
the classical propositional logic can take its beginning 
from the survival problem. It is also important that
this hypothetical origin of logic appears quite natural.

\section{About modal propositional logics} 

Following this idea, various types of logic can be 
viewed as theories of
such reductions for chosen classes of the operators. 
Here we suggest the following result concerning modal logic.
Its demonstration follows the expounded above method. 

{\sl{\bf Theorem 5.} Every propositional extention of the classical
propositional logic can be obtained by addition of one-place 
logical functions to the classical list $\to,\&,\vee,\neg$. 
}

{\sl\bf Proof.} Indeed, given $L=\{l_1,\dots,l_q\}$
and $M=L^n=L\times L\times\dots L$ we can construct
one-argument functions $\Gamma_i:L\to L$ where:
\begin{center}
$\Gamma_i(x) =
\begin{cases}
1,&  i\le x,\\
0,&   x<i.
\end{cases}$
\end{center}
Furthermore, we consider these functions $\Gamma_i,i=\overline{1,q},$ 
as functions
$\theta_i^j(x_1,\dots,x_j,\dots,\linebreak x_n),
i=\overline{1,q},j=\overline{1,n}$
such that $(\forall x_1\dots x_n)[\theta_i^j(x_1,\dots,x_j,\dots,x_n)=
\Gamma_i(x_j)]$. 
(Thus every $\theta_i^j$ has only one essential variable $x_j$.) These
functions $\theta_i^j$ satisfy the conditions of special 
$\theta$-functions above. Then we may take the closure relatively $\&$ of set 
$\{\theta_i^j|i=\overline{1,q},j=\overline{1,n}\}$ 
as the class $\cal K$ from theorem 3.
Hence we can use the theorem (as well as theorem 3$^{\star}$) 
to represent an arbitrary function
$\psi:M\to L$ by a formula constructed from standard operators $\to,\&,\vee$
and one-place functions $\Gamma_i,i=\overline{1,q}. \ \ \Box$

\section{Conclusions}

As the research demonstrates, the classical two-valued
propositional logic can be viewed merely as a realization of 
the above-mentioned principle of successive approximations for 
the decision-making problems within subject-environment survival framework. 

From this viewpoint, 
the classical propositional logic can take its beginning 
from the survival problem. It is also very important that
this hypothetical origin of logic appears quite natural. 

Then the approach can serve as a background for consideration of other 
families of mappings from one poset to another with a chosen notion of 
simplicity of mapping. Any such case 
generates a corresponding logic.

\end{document}